\input amstex
\magnification 1200
\documentstyle{amsppt}
\vcorrection{-1cm}

\topmatter
  \title
         Irreducibility of lemniscates
  \endtitle
  \author
         S.~Yu.~Orevkov
  \endauthor
\abstract
    We prove that lemniscates (i.e., sets of the form $|P(z)|=1$
    where $P$ is a complex polynomial)
    are irreducible real algebraic curves.
\endabstract
\endtopmatter

\def\refCC {1}
\def\refF  {2}

\def\R{\Bbb R}
\def\C{\Bbb C}
\def\P{\Bbb P}

\document
A {\it lemniscate} (or {\it polynomial lemniscate\/}) is a real
curve in $\C$ given by the equation $|P(z)|=1$ where $P(z)$ is
a non-constant polynomial with complex coefficients.

We say that a subset of $\R^2$ is an {\it irreducible real algebraic curve}
if it is the zero set of an irreducible over $\C$ real polynomial in two variables.
We prove in this note that

\medskip
\centerline{\sl any lemniscate is an irreducible real algebraic curve in $\C$}

\medskip
\noindent
(under the standard identification of $\R^2$ and $\C$).
This fact is an immediate consequence of Corollary 1 below.
(indeed, since $\{[P|=1\}=\{|P^d|=1\}$, any lemniscate
can be defined via a polynomal which is not a power of another one).

\proclaim{ Theorem 1 } Let $P$ and $Q$ be two polynomials in one variable
with complex coefficients. Then the polynomial $P(z)Q(w)-1$ is reducible
if and only if
$$
    P(z)=P_1(z)^d \qquad\text{and}\qquad Q(w)=Q_1(w)^d
$$
for $d>1$ and some polynomials $P_1(z)$ and $Q_1(w)$.
\endproclaim

\proclaim{ Corollary 1 }
Let $P(z)$ be a polynomial in one variable with complex coefficients and
 $f(x,y)$ be the real polynomial given by
$$
     f(x,y) = P(x+iy)\bar P(x-iy)-1
$$
where $\bar P$ is the polynomial whose coefficients are conjugates to those of $P$.
Then $f(x,y)$ is reducible over $\C$ if and only if $P(z)=P_1(z)^d$
for $d>1$ and a polynomial $P_1(z)$.
\endproclaim

\demo{ Proof } It is enough to apply Theorem 1 to $P$ and $\bar P$ after the linear
change of variables $z=x+iy$, $w=x-iy$ in $\C$.
\qed\enddemo

\medskip\noindent{\bf Remark 1.} (F.~B.~Pakovich).
If $P(z)$ and $Q(w)$ are arbitrary rational functions, then the problem of reducibility
of the algebraic curve $P(z)Q(w)=1$ seems to be very hard. For example, in another particular case (in a sense, opposite
to the ours), when $P(z)$ and $1/Q(w)$ are polynomials, this problem is solved in [\refF] (up to finite number of cases) and in [\refCC] (completely) but
the solution relies on the classification of simple finite groups.

\medskip
The rest of the paper is devoted to the proof of Theorem 1. Let
$$
  P(z)=\prod_{j=1}^k (z - z_j)^{p_j},
\;\; \deg P=p
\quad\text{and}\quad 
  Q(w)=\prod_{j=1}^l (w - w_j)^{q_j}, \;\; \deg Q=q
$$
where $z_1,\dots,z_k$ are pairwise distinct as well as $w_1,\dots,w_l$.
Let $C$ be the closure in $\P^1\times\P^1$ of the affine algebraic curve in $\C^2$
given by $P(z)Q(w)=1$ (here we represent $\P^1$ as $\C\cup\{\infty\}$).

Suppose that $P(z)Q(w)-1 = f'(z,w)f''(z,w)$ with non-constant polynomials $f'$ and $f''$.
Let $C'$ and $C''$ be the corresponding subsets of $C$.
We denote their local intersection numbers with the infinite lines
$L_1=\P^1\times\{\infty\}$ and
$L_2=\{\infty\}\times\P^1$ as follows:
$$
\xalignat 3
  &(C'\cdot L_1)_{(z_j,\infty)}=p'_j, &&(C''\cdot L_1)_{(z_j,\infty)}=p''_j, && (j=1,\dots,k),\\
  &(C'\cdot L_2)_{(\infty,w_j)}=q'_j, &&(C''\cdot L_2)_{(\infty,w_j)}=q''_j, && (j=1,\dots,l).
\endxalignat
$$
Let $(p',q')$ and $(p'',q'')$ be the bidegree of $C'$ and $C''$ respectively.
Then
$$
    p=p'+p'',\quad q=q'+q'',\quad p_j=p'_j+p''_j,\quad q_j=q'_j+q''_j.
$$

The germ of $C$ at $(z_j,\infty)$ has equation $u^q=v^{p_j}$
in some local analytic coordinates $(u,v)$. So, it has $\gcd(q,p_j)$ local branches
which are distributed in some proportion between $C'$ and $C''$.
By comparing the the degree of the projections
of the germs of $C'$ and $C''$ onto the coordinate axes, we conclude that
$p'_j/p''_j = q'/q''$. Similarly, $q'_j/q''_j= p'/p''$.
We denote these quotients by $\alpha$ and $\beta$ respectively.
Since $\sum p'_j=p'$ and $\sum p''_j=p''$, we obtain
$$
   \beta p'' = p' = p'_1+\dots+p'_k = \alpha p''_1+\dots+\alpha p''_k = \alpha p''
$$
whence $\alpha=\beta$. Let $\alpha=d'/d''$ with coprime $d'$ and $d''$.
Since
$$
    \frac{p'_1}{p''_1}=\dots=\frac{p'_k}{p''_k}=
    \frac{q'_1}{q''_1}=\dots=\frac{q'_l}{q''_l}=\frac{d'}{d''},
$$
we obtain $p'_j=a_j d'$, $p''_j=a_j d''$ and $q'_j=b_j d'$, $q''_j=b_j d''$
for some integers $a_1,\dots,a_k$ and $b_1,\dots,b_l$.
Hence $p_j = p'_j+p''_j=a_j d$ and $q_j=b_j d$ for $d=d'+d''$,
and we finally obtain $P(z)=P_1(z)^d$ and $Q(w)=Q_1(w)^d$ with
$$
    P_1(z)=\prod_{j=1}^k (z - z_j)^{a_j}
    \quad\text{and}\quad 
    Q_1(w)=\prod_{j=1}^l (w - w_j)^{b_j}.
$$

I am gateful to Fedor Pakovich for suggesting the problem and
stimulating discussions.

\enddocument